\documentclass[a4paper]{article}
\usepackage{amsmath}
\usepackage{lmodern}
\usepackage{amsmath}
\usepackage{amssymb}
\usepackage{mathrsfs}
\usepackage{graphicx}

\newtheorem{theorem}{Theorem}

\newtheorem{corollary}[theorem]{Corollary}

\newtheorem{definition}[theorem]{Definition}
\newtheorem{example}[theorem]{Example}

\newtheorem{lemma}[theorem]{Lemma}

\newtheorem{remark}[theorem]{Remark}

\newcommand{\proof}{\noindent\textbf{Proof.}\quad }
\newcommand{\R}{{\mathbb R}}
\newcommand{\N}{{\mathbb N}}
\def\qed{\hbox to 0pt{}\hfill$\rlap{$\sqcap$}\sqcup$}

\begin{document}
	\title{Nonlinear Differential Equations with Perturbed Dirichlet Integral Boundary Conditions$^1$}
	\author{Alberto Cabada and Javier Iglesias\\
		$^1$Departamento de Estat\'istica, An\'alise Matem\'atica e Optimizaci\'on\\
		Instituto de Matem\'aticas, Facultade de Matem\'aticas, \\
		Universidade de Santiago de Compostela, Spain.\\
		alberto.cabada@usc.gal, \quad iglesiasprez.javier@gmail.com }
	\maketitle
	\begin{abstract}
		This paper is devoted to prove the existence of positive solutions of a second order differential equation with a  nonhomogeneous Dirichlet conditions given by a parameter dependence integral. The studied problem is a nonlocal perturbation of the Dirichlet conditions by considering a homogeneous Dirichlet-type condition at one extreme of the interval and an integral operator on the other one. We obtain the expression of the Green's function related to the linear part of the equation and characterize its constant sign. Such property will be fundamental to deduce the existence of solutions of the nonlinear problem. The results hold from fixed point theory applied to related operators defined on suitable cones.
	\end{abstract}
	
	
	
	\section{Introduction}

This paper is devoted to the study of the existence of solutions of the following family of nonlinear second order ordinary differential equations
	\begin{equation}\label{ecuacion_general}
		u''(t)+\gamma u(t)+f(t,u(t))=0,\qquad 0<t<1,
	\end{equation}
	coupled to the following integral boundary conditions
	\begin{equation}\label{condiciones}
		u(0)=0,\qquad u(1)=\lambda\int_{0}^{1}u(s)ds,
	\end{equation}
	where $ \gamma <\pi ^ 2 $ and $ \lambda \geq 0 $.

	To this end, we will distinguish the cases $ \gamma = 0 $, $ \gamma> 0 $ and $ \gamma <0 $. We will analyze each of them and give optimal sufficient conditions on $\gamma$, $ \lambda $ and $ f $ that allow us to ensure the existence of a solution of the considered problem.

This kind of problems model the behavior of an harmonic oscillator, subject to a external force $f$, which is fixed at the left extreme of the interval and has some mechanism at the right one, that controls the displacement according to the feedback from devices measuring the displacements along parts of the oscillator. Integral boundary conditions have been considered in many works in the literature, see for instance, \cite{Cabada-Jebari, Hu-Yan, Zhang-Abella-Feng} (for second and fourth order Ordinary Differential Equations) or \cite{Ahmad-Hamdan-Alsaedi-Ntouyas, Cabada_articulo,  Chandran-Gopalan-Tasneem-Abdeljawad, Duraisamy-Nandha-Subramanian} (for Fractional equations) and references therein.

We refer a function $u\in \mathcal {C} ([0,1])$ as a nonnegative solution of problem \eqref{ecuacion_general}-\eqref{condiciones} if $u$ solves such problem and $u(t)\geq 0$, for all $t\in [0,1]$. A function $u\in \mathcal {C} ([0,1])$ is called positive solution of problem \eqref{ecuacion_general}-\eqref{condiciones} if it is a nonnegative solution and $u(t)>0$, for all $t\in(0,1)$.

The paper is organized as follows: in Section 2, we study the linear part of problem \eqref{ecuacion_general}--\eqref{condiciones}, where we obtain the explicit expression of the related Green's function and calculate the exact values of $\gamma$ and $\lambda$ for which the Green's function has constant sign. In next section we prove the existence of positive solutions for the nonlinear problem \eqref{ecuacion_general}--\eqref{condiciones}. Such solutions are given as the fixed points of a related integral operator defined on a suitable cone. At the end of this section we show two examples where the applicability of the obtained results is pointed out.

The following concept will be fundamental in order to deduce our existence results.

\begin {definition}
Let $ X $ be a Banach space. A subset $ K \subset X $ is a cone if:
\begin {itemize}
\item $ K $ is closed
\item $ K + K \subset K, \hspace {1 ex} \lambda K \subset K $ for all $ \lambda \ge0 $ and $ K \cap (-K) = \{0 \} $.
	\end {itemize}
			\end {definition}

		We will use the celebrated expansion/contraction theorem of Krasnosels'ki\u{\i} \cite{Kras_general}:
		\begin {theorem} [Krasnosels'ki\u{\i}] \label{Krasnoselskii}
		Let $ X $ be a Banach space and $ K \subset X $ a cone in $ X $. Let $ \Omega_1, \Omega_2 \subset X $ open bounded such that $ 0 \in \Omega_1 \subset \overline {\Omega_1} \subset \Omega_2 $ and
		$ T: K \cap (\overline {\Omega_2} \setminus \Omega_1) \to K $ a compact operator that satisfies one of the following properties:
		\begin {enumerate}
		\item $ || T (u) || \geq || u ||, \hspace {1 ex} \forall u \in K \cap \partial \Omega_1 $ y $ || T (u) || \leq || u ||, \hspace {1 ex} \forall u \in K \cap \partial \Omega_2 $.
		\item $ || T (u) || \leq || u ||, \hspace {1 ex} \forall u \in K \cap \partial \Omega_1 $ y $ || T (u) || \geq || u ||, \hspace {1 ex} \forall u \in K \cap \partial \Omega_2 $.
		\end {enumerate}
		Then $ T $ has a fixed point at $ K \cap (\overline {\Omega_2} \setminus \Omega_1) $.
		\end {theorem}
%
%
%

\section{Linear part: Green's function}

In this section we obtain the expression of the related Green's function of the linear part of problem \eqref{ecuacion_general}--\eqref{condiciones} and deduce some important properties that will be fundamental to obtain the existence of positive solutions of the nonlinear problem. To this end, we consider separately three cases depending on the sign of the real parameter $\gamma$.  To do this, we will follow a treatment of a similar problem studied in \cite{Cabada_articulo} for fractional equations.

	\subsection{Case $\gamma=0$}\label{Secc:gamma=0}

In this subsection, we obtain the expression of the Green's function related to the linear problem	
\begin{equation}\label{ecuacion1_lineal}
		u''(t)+\sigma(t)=0, \quad 0<t<1,
	\end{equation}
coupled to the boundary conditions (\ref{condiciones}).
	
	\begin{theorem}\label{ecuacion1_green_teorema}
Let $ \lambda \neq 2 $ and $ \sigma \in \mathcal{C}( [0,1]) $, then  problem (\ref{ecuacion1_lineal}),  (\ref{condiciones}) has a unique solution $ u \in \mathcal {C}^2([0,1]) $, which is given by the following expression
		
		$$u(t)=\int_{0}^{1}G(t,s)\sigma(s)ds,$$
	where
		\begin{equation}\label{ecuacion1_green}
			G(t,s)=\begin{cases}
				\frac{t(1-s)(2-\lambda+\lambda s)-(2-\lambda)(t-s)}{2-\lambda},
				& 0\le s\le t\le 1,\\
				\frac{t(1-s)(2-\lambda+\lambda s)}{2-\lambda}, 
				& 0\le t<s\le 1.
			\end{cases}
		\end{equation} 
			\end{theorem}
		
		\proof
By using the Fundamental Theorem of Integral Calculus to the equation (\ref{ecuacion1_lineal}), together with Fubini's Theorem, we arrive at the following expression	
\[
u(t)
=-\int_{0}^{t}(t-s)\sigma(s)ds+c_1t+c_2.
\]	

			Since $u(0)=0$, we deduce that $c_2=0$.
			
			Now, the boundary condition at $t=1$ implies that
			$$\lambda\int_{0}^{1}u(s)ds=u(1)=-\int_{0}^{1}(1-s)\sigma(s)ds+c_1$$
			so
			$$c_1=\int_{0}^{1}(1-s)\sigma(s)ds+\lambda\int_{0}^{1}u(s)ds.$$
		
		As a consequence 
			\begin{equation}\label{ecuacion1_solucion}
				u(t)=-\int_{0}^{t}(t-s)\sigma(s)ds+t\int_{0}^{1}(1-s)\sigma(s)ds+\lambda t\int_{0}^{1}u(s)ds.
			\end{equation}
			By denoting $A=\int_{0}^{1}u(s)ds$, we have that 
			$$A=\int_{0}^{1}u(t)dt=-\int_{0}^{1}\int_{0}^{t}(t-s)\sigma(s)dsdt+\int_{0}^{1}\int_{0}^{1}t(1-s)\sigma(s)dsdt+\lambda A\int_{0}^{1}tdt.$$
		
	Thus, 
			$$\begin{aligned}
				A=&-\int_{0}^{1}\int_{s}^{1}(t-s)\sigma(s)dtds+\int_{0}^{1}\int_{0}^{1}t(1-s)\sigma(s)dtds+\frac{\lambda}{2}A\\
				=&-\int_{0}^{1}\frac{(1-s)^2}{2}\sigma(s)ds+\int_{0}^{1}\frac{(1-s)}{2}\sigma(s)ds+\frac{\lambda}{2}A.\end{aligned}$$
			
			Then
			
			 $$\begin{aligned}
				A=&\frac{2}{2-\lambda}\left(-\int_{0}^{1}\frac{(1-s)^2}{2}\sigma(s)ds+\int_{0}^{1}\frac{(1-s)}{2}\sigma(s)ds\right)\\
				=&-\frac{1}{2-\lambda}\int_{0}^{1}(1-s)^2\sigma(s)ds+\frac{1}{2-\lambda}\int_{0}^{1}(1-s)\sigma(s)ds.
			\end{aligned}$$
		
Substituting this value in (\ref{ecuacion1_solucion}), we obtain the expression of the solution $u$ as follows:
			$$\begin{aligned}
				u(t)=&-\int_{0}^{t}(t-s)\sigma(s)ds+t\int_{0}^{1}(1-s)\sigma(s)ds-\frac{\lambda}{2-\lambda}t\int_{0}^{1}(1-s)^2\sigma(s)ds\\
				&+\frac{\lambda}{2-\lambda}t\int_{0}^{1}(1-s)\sigma(s)ds\\
				=& -\int_{0}^{t}(t-s)\sigma(s)ds+t\int_{0}^{1}\frac{(1-s)(2+
					\lambda(s-1))}{2-\lambda}\sigma(s)ds\\
				=&\int_{0}^{t}\frac{t(1-s)(2+\lambda s-\lambda)) -(2-\lambda)(t-s)}{2-\lambda} \sigma(s)ds\\
				&+\int_{t}^{1}\frac{t(1-s)(2+\lambda s-\lambda))}{2-\lambda}\sigma(s)ds=\int_{0}^{1}G(t,s)\sigma(s)ds.
			\end{aligned}$$
		\qed

	
In the sequel, we will state two lemmas related to the properties of the Green's function that will be useful to prove the existence of a positive solution of the nonlinear  problem (\ref{ecuacion_general})--(\ref{condiciones}) with $\gamma =0$.
	\begin{lemma}\label{prop_green1}
		Let $G$ be the Green's function related to  problem (\ref{ecuacion1_lineal}),  (\ref{condiciones}), given by expression (\ref{ecuacion1_green}). Then, for all   $\lambda\neq 2$ the following properties are fulfilled:
		\begin{enumerate}
			\item $G(0,s)=G(t,0)=G(t,1)=0$,  for all  $t,s\in[0,1].$
					\item $G(t,s)$ is continuous on $[0,1]\times[0,1].$	\item $G(1,s)=0$  for all  $s\in(0,1)$  if and only if  $\lambda=0.$
			\item $(2-\lambda)G(1,s)>0$  for all  $s\in(0,1)$ if and only if $\lambda > 0$.
			\item $(2-\lambda)G(s,s)>0$  for all  $s\in(0,1)$ if and only if $\lambda < 2$.
			\item $G(t, s)> 0 $ for all $ t, s \in (0,1)$ if and only if $ \lambda  \in \left[0,2\right)$.
			\item $G(t, s) $ changes sign on $(0,1 ) \times (0,1)$  for all $ \lambda \not \in \left[0,2\right)$.
			\item For all  $\lambda\in [0,2),\hspace{1ex} G(t,s)\leq\frac{1}{2(2-\lambda)},\quad\forall t,s\in[0,1].$	
		\end{enumerate}
		\end{lemma}
	
		\proof
	Properties 1. and 2. are immediate. Let's now prove the remaining properties:
			\begin{enumerate}
				\item[3.] Let $s\in(0,1)$, then $G(1,s)=0$ if and only if $ \lambda s(1-s)=0$, i. e. $ \lambda=0$.
				
				\item[4.] The result holds trivially from the fact that 
				\begin{equation}
					\label{e-G(1,s)}
						(2-\lambda)G(1,s)=\lambda s(1-s), \hspace{1ex}\text{ for all  }s\in(0,1).
				\end{equation}
			
				\item[5.] The result is immediately deduced from the following equality:
					\begin{equation}
					\label{e-G(s,s)}
					 (2-\lambda)G(s,s)=s\,(1-s)\, (2-\lambda\,(1-s)) , \hspace{1ex}\text{ for all  }s\in(0,1).
					 \end{equation}
					\item[6.] Since $G(t,s)$ is linear on $t$, for all $s\in[0,1]$ fixed, $G(t,s)$ attains its maximum and minimum at $t=0,t=s$ or at $t=1$.\\
				From Property 1., we have that $G(0,s)=0$  for all  $s\in[0,1].$\\
				From Property 5., we have that $G(s,s)>0$ for all  $s\in (0,1)$, if and only if $\lambda <2$.\\
				From Property 4., we have that $G(1,s)>0$ for all  $s\in (0,1)$, if and only if $\lambda \in (0,2)$.\\
				As a consequence of the three previous assertions and Property 3., this Property holds.
				
				\item[7.]
				
				From \eqref{e-G(1,s)} it is clear that $G(1,s)<0$ for all  $s\in (0,1)$, if and only if $\lambda \not \in (0,2)$.\\
				From \eqref{e-G(s,s)} it is clear that $G(s,s)>0$ for all  $s\in (0,1)$ and $\lambda \le 0$ and, moreover, if $\lambda >2$, we have that $G(s,s)>0$ for $s\in (0,1)$ close enough to $0$. Thus, this property is fulfilled.
				
					\item[8.] From Property 6., we know that $G(t,s)>0$ for all $(t,s) \in (0,1) \times (0,1)$. As in Property 6., we know that the maximum values will be attained at $G(s,s)$ and/or $G(1,s)$. Now, since
				$$(2-\lambda)G(s,s)=s(1-s)(2-\lambda(1-s))=2s(1-s)-\lambda s(1-s)^2\leq 2s(1-s)\leq\frac{1}{2}$$ 
			and
				$$(2-\lambda)G(1,s)=\lambda s(1-s)\leq \frac{\lambda}{4}<\frac{1}{2},$$
				the proof is concluded.	
			\end{enumerate}  
		\qed

In the sequel, we deduce two sharp inequalities for the positiveness of the Green's function.

	\begin{lemma}\label{cotas_green1_enunc}
		Let $\lambda\in(0,2)$ and $G(t,s)$ be the Green's function related to problem (\ref{ecuacion1_lineal}),  (\ref{condiciones}), given by expression (\ref{ecuacion1_green}). Then the following properties hold:
		\begin{equation}\label{cotas_green1}
			t\,G(1,s)\leq G(t,s)\leq \frac{2}{\lambda}G(1,s), \quad\text{ for all  } s,t\in[0,1].
		\end{equation}
			\end{lemma} 
		\proof
		For $ t = 0$ $t=1$, $s = 0$ or $ s = 1 $ the inequalities follow immediately from  Properties 1. and 4. in  Lemma \ref{prop_green1}. 
		
		Now let $ t, s $ be such that $ 0 <t \leq s <1 $. In this case:	
		$$h(t,s):=\frac{G(t,s)}{G(1,s)}=\frac{t(2-\lambda(1-s))}{\lambda s}=t \left(1+\frac{2-\lambda}{\lambda s}\right).$$
		As a consequence, we have
			$$t<t\left(1+\frac{2-\lambda}{\lambda s}\right)= h(t,s)\leq t\frac{2}{\lambda s}=\frac{2}{\lambda}\frac{t}{s}\leq \frac{2}{\lambda}.$$
			
			Consider now the case $0<s\le t<1$, using that $s\geq ts$, we conclude that
			$$\begin{aligned}
				h(t,s)=&\frac{t(1-s)(2-\lambda(1-s))-(2-\lambda)(t-s)}{\lambda s(1-s)}\\
				\geq&\frac{t(1-s)(2-\lambda(1-s))-(2-\lambda)(t-ts)}{\lambda s(1-s)}\\
				=& \frac{t(1-s)[2-\lambda(1-s)-(2-\lambda)]}{\lambda s(1-s)}=t.
			\end{aligned}$$
		
		Moreover,
			$$\lim\limits_{t\to s^+}h(t,s)= \frac{2-\lambda(1-s)}{\lambda}\leq\frac{2}{\lambda}$$
			and
			$$ \lim\limits_{t\to 1^-}h(t,s)= \frac{(1-s)[(2-\lambda(1-s))-(2-\lambda)]}{\lambda s(1-s)}=1.$$
			
			On the other hand, since $h(t,s)$ is a non-negative, continuous and linear function with respect to $t$, we deduce 
			$$t \le  h(t,s)\leq \max{\left\{\lim\limits_{t\to 1^-}h(t,s),\lim\limits_{t\to s^+}h(t,s)\right\}}=\frac{2}{\lambda}.$$
			From Property $4.$ in Lemma \ref{prop_green1}, we conclude that the inequalities (\ref{cotas_green1}) are satisfied.
		\qed

	\subsection{Case $\gamma>0$}

In this subsection, we will obtain the expression of the Green's function related to the problem
\begin{equation} \label{equation2_linear}
	u '' (t) + m ^ 2u (t) + \sigma (t) = 0, \quad 0 <t <1,
\end{equation}
coupled to the boundary conditions (\ref{condiciones}).

\begin{theorem} \label{equation2_green_teorem}
	Let $ \lambda \neq \frac {m \sin m} {1- \cos m} $, $m>0$, $m \neq 2 \, k\, \pi $, $k=1, 2, \ldots$, and $ \sigma \in \mathcal {C} ([0,1]) $. Then  problem (\ref{equation2_linear}),  (\ref{condiciones}) has a unique solution $ u \in \mathcal {C} ^ 2 ([0,1]), $ which is given by $$ u (t) = \int_{0} ^ {1} G_m (t, s) \sigma (s) ds, $$ where
	\begin{equation} \label{equation2_green}
		G_m (t, s) = \begin{cases}
			G ^ 1_m (t, s), & 0 \le s \le t \le 1, \\
			G ^ 2_m (t, s), & 0 \le t <s \le 1.
		\end{cases}
	\end{equation}
	Here, if $ m \neq k\, \pi$, $k \in \N$ odd, 
	$$ \begin{aligned}
		G ^ 1_m (t, s) = & \frac {\sin (ms) [\sin (m-mt) (m \sin m- \lambda (1- \cos m)) + \lambda \sin (mt) ]} {m \sin m (m \sin m- \lambda (1- \cos m))} \\
		& + \frac {\lambda \sin (mt) (\sin (m-ms) - \sin m)} {m \sin m (m \sin m- \lambda (1- \cos m))}
	\end{aligned} $$ and $$
	G ^ 2_m (t, s) = \frac {\sin (mt) [\sin (m-ms) (m \sin m + \lambda \cos m) + \lambda (\sin (ms) - \sin m) ]} {m \sin m (m \sin m- \lambda (1- \cos m))}, $$
	and, if $ m = k\, \pi$, for some $k \in \N$ odd,
		$$ \begin{aligned}
		G ^ 1_{ k\, \pi} (t, s) = & \frac{2 \lambda  \sin (\pi  k s) \cos (\pi  k t)+\sin (\pi  k t) (\lambda  (-\cos (\pi  k s))-\pi  k \sin (\pi  k s)+\lambda )}{2 \pi  k \lambda }
	\end{aligned} $$ and $$
	G ^ 2_{ k\, \pi}  (t, s) = \frac{\sin (\pi  k t) (\lambda  \cos (\pi  k s)-\pi  k \sin (\pi  k s)+\lambda )}{2 \pi  k \lambda }.$$

\end{theorem}
\proof
It is immediate to verify that the spectrum of problem (\ref{equation2_linear}),  (\ref{condiciones}) is given by the following pairs on the plane $(m,\lambda)$:

\begin{enumerate}
	\item $\left(m,\frac {m \sin m} {1- \cos m}\right), \quad m> 0, \; m \neq 2 \, k\, \pi, \; k=1, 2, \ldots$
	\item $\left(2 \, k\, \pi, \lambda\right), \; \lambda \in \R,\;k=1, 2, \ldots$
\end{enumerate}
Consider the case  $m \neq k\, \pi$, $k \in \N$ odd, and let $ v $ be the unique solution of
$$ \begin{cases}
	v '' (t) + m ^ 2v (t) + \sigma (t) = 0, \quad 0 <t <1, \\
	v (0) = 0, v (1) = 0,
\end{cases} $$
and $ w $ be defined as the unique solution of
$$ \begin{cases}
	w '' (t) + m ^ 2w (t) = 0, \quad 0 <t <1, \\
	w (0) = 0, w (1) = 1,
\end{cases} $$
then, it is not difficult to verify that  $ u (t) = v (t) + \left(\lambda \int_{0} ^ {1} u (s) ds\right) w (t) $ is a solution of (\ref{equation2_linear}),  (\ref{condiciones}).\\

Using \cite {programa_green} we have that
$$ v (t) = \int_{0} ^ {1} G ^ v_m (t, s) \sigma (s) ds, $$ where
$$ G ^ v_m (t, s) = \begin{cases}
	\frac {\sin (ms) \sin (m-mt)} {m \sin m}, & 0 \le s \le t \le1, \\\frac {\sin (m-ms) \sin (mt )} {m \sin m}, & 0 \le t <s \le 1.
\end{cases} $$

It is immediate to verify that $$ w (t) = \frac {\sin (mt)} {\sin m} . $$

Thus, \begin{equation} \label{equation2_solution}
	u (t) = \int_{0} ^ {1} G ^ v_m (t, s) \sigma (s) ds + \frac {\lambda} {\sin m} \left (\int_{0} ^ {1 } u (s) ds \right) \sin (mt).
\end{equation}

Denoting $ A = \int_{0} ^ {1} u (s) ds $, we deduce from the previous expression that
$$ \begin{aligned}
	A = & \int_{0} ^ {1} u (t) dt = \int_{0} ^ {1} \int_{0} ^ {1} G ^ v_m (t, s) \sigma (s) dsdt + \frac {\lambda} {\sin m} A \int_{0} ^ {1} \sin (mt) dt \\
	= & \int_{0} ^ {1} \int_{0} ^ {t} \frac {\sin (ms) \sin (m-mt)} {m \sin m} \sigma (s) dsdt + \int_{0} ^ {1} \int_{t} ^ {1} \frac {\sin (m-ms) \sin (mt)} {m \sin m} \sigma (s) dsdt \\& + \frac {\lambda} {\sin m} A \int_{0} ^ {1} \sin (mt) dt = \int_{0} ^ {1} \frac {\sin (ms) (1- \cos ( m-ms))} {m ^ 2 \sin m} \sigma (s) ds \\& + \int_{0} ^ {1} \frac {\sin (m-ms) (1- \cos (ms ))} {m ^ 2 \sin m} \sigma (s) ds + A \frac {\lambda} {\sin m} \frac {1- \cos m} {m}. \end{aligned} $$ 

As a consequence, it follows that 
$$ \begin{aligned}
	A & = \left (1- \frac {\lambda (1- \cos m)} {m \sin m} \right) ^ {- 1} \int_{0} ^ {1} \frac {\sin (ms ) + \sin (m-ms) - \sin m} {m ^ 2 \sin m} \sigma (s) ds \\
	& = \int_{0} ^ {1} \frac {\sin (ms) + \sin (m-ms) - \sin m} {m (m \sin m- \lambda (1- \cos m)) } \sigma (s) ds.
\end{aligned} $$

Substituting this value in (\ref{equation2_solution}) we obtain the following expression
$$ \begin{aligned}
	u (t) = & \int_{0} ^ {t} \frac {\sin (ms) \sin (m-mt)} {m \sin m} \sigma (s) ds + \int_{t} ^ { 1} \frac {\sin (m-ms) \sin (mt)} {m \sin m} \sigma (s) ds \\& + \int_{0} ^ {1} \frac {\sin (ms ) + \sin (m-ms) - \sin m} {m (m \sin m- \lambda (1- \cos m))} \frac {\lambda \sin (mt)} {\sin m} \sigma (s) ds \\
	= & \int_{0} ^ {t} \left (\frac {\sin (ms) [\sin (m-mt) (m \sin m- \lambda (1- \cos m)) + \lambda \sin (mt)]} {m \sin m (m \sin m- \lambda (1- \cos m))} \right. \\& \left. + \frac {\lambda \sin (mt) (\sin (m-ms) - \sin m)} {m \sin m (m \sin m- \lambda (1- \cos m))} \right) \sigma (s) ds \\
	& + \int_{t} ^ {1} \frac {\sin (mt) [\sin (m-ms) (m \sin m + \lambda \cos m) + \lambda (\sin (ms) - \sin m)]} {m \sin m (m \sin m- \lambda (1- \cos m))} \sigma (s) ds \\
	= & \int_{0} ^ {1} G_m (t, s) \sigma (s) ds
\end{aligned} $$
The uniqueness of the Green's function is deduced from the uniqueness of functions $ v $ and $ w $.

The expressions of the Green's function for $m=k\, \pi$, with $k \in \N$ odd, follow by tacking the limit of the expressions of $G_m^1$ and $G_m^2$ when $m$ goes to $k\, \pi$. By direct calculations, it is immediate to verify that such function satisfies the properties of the Green's function of problem (\ref{equation2_linear}),  (\ref{condiciones}) with $m=k\, \pi$,  $k \in \N$ odd.
\qed
\\

In the same way as in the case $ \gamma = 0 $, we will now state different results about the properties of the Green's function that we have just obtained. As before, we will try to give conditions that allow us to ensure that the Green's function attains a constant sign. To this end, we will use the following property, which is a direct consequence of the Sturm comparison theorem


\begin{remark} \label{distancia_ceros}
	Since $ w (t) = \sin (\pi t), $ is a nontrivial solution of $$ w '' (t) + \pi ^ 2w (t) = 0, $$ and it vanishes in $ \mathbb {Z} $, from the classical Sturm comparison result, we have that for $ \gamma <\pi ^ 2 $ any non-trivial solution of the problem $$ v '' (t) + \gamma\, v (t) = 0 $$ vanishes at most once in the interval [0,1].	
	Therefore, the distance between two zeros of any non-trivial solution must be greater than $ 1 $.
\end{remark}

\begin{lemma} \label{prop_green2}
	Let $ G_m $ be the Green's function related to  problem (\ref{equation2_linear}),  (\ref{condiciones}), given by the expression (\ref{equation2_green}). Then for all $ \lambda \neq \frac {m \sin m} {1- \cos m} $, $m>0$, $m \neq 2 \, k\, \pi $, $k=1, 2, \ldots$, the following properties hold:
	\begin{enumerate}
		\item $ G_m (0, s) = G_m (t, 1) =  G_m (t, 0)=0$, for all $ t, s \in [0,1]. $
		\item $ G_m (t, s) $ is continuous at $(t, s) \in  [0,1] \times [0,1]. $
		\item If $m \in (0,2\, \pi)$, then $ G_m(1, s) = 0 $ for all $ s \in (0,1) $ if and only if $ \lambda = 0 $.
		\item If $m \in (0,2\, \pi)$, then $ (m \sin m- \lambda (1- \cos m)) G_m (1, s)> 0 $ for all $ s \in (0,1)$ and $\lambda>0. $
		\item $G_m (t, s)> 0 $ for all $ t, s \in (0,1 )$ if and only if $ 0 \le \lambda <\frac {m \sin m} {1- \cos m}$, $m \in (0, \pi]$.
		\item $G_m (t, s) $ changes sign on $(0,1 ) \times (0,1)$  for all $ \lambda \not \in \left[0,\frac {m \sin m} {1- \cos m}\right)$, $m \in (0, \pi]$.
	\end{enumerate}
\end{lemma}
\proof
Properties 1. and 2. are immediate. Let's now see the others:
\begin{enumerate}
	\item [3.] Let $ s \in (0,1) $, then $G_m (1, s) = 0$ if and only if
	$$
	\frac {\lambda (\sin (ms) + \sin (m-ms) - \sin m)} {m (m \sin m- \lambda (1- \cos m))} = 0
	$$
	which is equivalent to
	$$ \lambda (\sin (ms) + \sin (m-ms) - \sin m) =: \lambda\, r_m(s)= 0.$$
	
	It is easy to check that if $m \in (0,2\, \pi)$ then function $ r_m$
	has a unique maximum at $ s = \frac {1} {2} $ and, since $ r_m (0) = 0 = r_m (1) $, we deduce that $ r_m (s)> 0 $ for $ s \in (0,1) $. Therefore $ G_m (1, s) = 0 $ if and only if $ \lambda = 0 $.
	\item [4.]  From previous assertion, we have that 
	$$ (m \sin m- \lambda (1- \cos m)) G_m (1, s) = \frac {\lambda} {m} r_m(s)> 0,$$
	for all $ s \in (0,1)$, $m \in (0,2\, \pi)$ and $\lambda >0$.
	\item [5.] Assume that $ 0 \le \lambda <\frac {m \sin m} {1- \cos m}$, $m \in (0, \pi]$. Using properties 1. and 4., we know that for each $ s \in (0,1) $, $ G_m (0, s) = 0 $ and $ G_m (1, s)> 0 . $ In addition,
	$$ \begin{aligned}
		\frac {\partial G_m} {\partial t} (0, s) = & \frac {\sin (m(1-s)) } {\sin m }  + \frac {\lambda \,r_m(s)} {\sin m ( m \sin m + \lambda (1- \cos m))}> 0.
	\end{aligned} $$
	Therefore, $ G_m (t, s) $ is increasing and positive to the right of $ t = 0 $ for every $ s \in (0,1) $. \\
	On the other hand, 
	$$ \frac {\partial ^ 2 G_m} {\partial t ^ 2} (t, s) + m ^ 2G_m (t, s) = 0,\; t \in [0,1], \quad t \neq s. $$ 
	
	Since $ G_m (t, s) \in \mathcal {C} ^ 2 (\left [0 , s \right) \cup \left (s, 1 \right]) $, we can use Remark \ref{distancia_ceros} on each of the intervals. \\
	Suppose that there exists $ t_0 \in (0,1) $ such that $ G_m (t_0, s) = 0 $. We distinguish the following cases:
	\begin{itemize}
		\item If $ t_0 \in \left [0, s \right) $, since $ G_m (0, s) = 0 $, we would have two zeros that distances less than 1, which is not possible.
		\item Now suppose $ t_0 = s $, in this case we can continuously extend $ G_m (t, s) $ to the interval $ [0, s] $, which brings us back to the previous case.
		\item If $ t_0 \in (s, 1) $, since $ G_m (s, s)> 0 $ and $G_m(1,s)\ge 0$, can occur three situations: \begin{enumerate}
			\item If $ \frac {\partial G_m} {\partial t} (t_0, s) <0 $, then there exists $ t_1 \in ( t_0,1] $ such that $ G_m (t_1, s) = 0 $, which is not possible as a consequence of Remark \ref{distancia_ceros}.
			\item If $ \frac {\partial G_m} {\partial t} (t_0, s) = 0 $, then $ G_m (t, s) \equiv 0 $ in $ \left (s, 1 \right] $, which is not possible.
			\item If $ \frac {\partial G_m} {\partial t} (t_0, s)> 0 $, then there exists $ t_1 \in (s, t_0) $ such that $ G (t_1, s) = 0 $, reaching a contradiction again.
		\end{enumerate}
		From the equality  
		$$ (m \sin m- \lambda (1- \cos m)) G_m (1, s) = \frac {\lambda} {m} r_m(s),$$
		we deduce immediately that $ G_m (1, s)$ takes positive and negative values for all $m >2\, \pi$, $m \neq 2\, k\, \pi$, $k=1,2, \ldots$ and all $\lambda \in \R$, $\lambda \neq 0$.
		Moreover,  from the expression of $G_m(1,s)$ we deduce that if either, $\lambda<0$ and $m \in (0,2\, \pi)$, or $ \lambda >\frac {m \sin m} {1- \cos m}$ and $m \in (0,
		\pi]$, or $\lambda >0$ and $m \in (\pi,2\, \pi)$, then $G_m(1,s)<0$.
		\item [6.] Using previous assertion, we only need to verify that if  $\lambda \not \in \left[0,\frac {m \sin m} {1- \cos m}\right)$, $m \in (0, \pi]$, then $G_m(t,s)$ takes some positive values on $(0,1) \times (0,1)$. But, to verify this property it is enough to consider function 
		\[
		f_m(t):= G_m(t,t).
		\]
		
		By direct computation, we have that $f'_m(0)=1$ for all $\lambda \in \R$ and $m >0$, $m \neq 2\, k\, \pi$, $k=1,2, \ldots$.  In consequence, $G_m(t,t) >0$ in a small enough neighborhood of $(0,0)$.
	\end{itemize}
\end{enumerate}
\qed

Now we deduce the following stronger condition on the Green's function

\begin{lemma}\label{cotas_green2_enunc}
	Let $0<m<\pi,\hspace{1ex}0\le\lambda<\frac{m\sin m}{1-\cos m}$ and $G_m(t,s)$ be the Green's function of problem (\ref{equation2_solution}), (\ref{condiciones}) given by expression (\ref{equation2_green}). Then, there are $h_1\in\mathcal{C}([0,1]),\hspace{1ex}h_1>0$ on $\left(0,1\right]$ and $C_1\in\mathbb{R},\hspace{1ex}C_1>0$ such that:
	\begin{equation} \label{cotas_green2}
		h_1 (t) G_m (t, s) \leq G_m (t, s) \leq C_1\,G_m (1, s), \qquad \text {for all} \hspace {1ex} t, s \in [0,1]
	\end{equation}
\end{lemma}
\proof
If $ t = 0, \hspace {1ex} s = 0 $ or $ s = 1 $ the result follows from Lemma \ref{prop_green2}. Let then $ 0 <t \leq1 $ arbitrarily set.
Since $$ \lim \limits_ {s \to 0 ^ +} \frac {G_m (t, s)} {G_m (1, s)} = \frac {\sin (m-mt) (m \sin m - \lambda (1- \cos m))} {\lambda \sin m (1- \cos m)}> 0, $$
$$ \lim \limits_ {s \to 1 ^ -} \frac {G_m (t, s)} {G_m (1, s)} = 1 $$ and taking into account the properties 4. and 5. of  Lemma \ref{prop_green2}, then, for any $s \in (0,1)$ fixed, we can extend function $ \frac {G_m (\cdot, s)} {G_m (1, s)} $ continuously to the interval $ [0,1] $ and, furthermore, this extension is strictly positive for all $ t \in \left (0,1 \right] $. \\
As a consequence, there are $$ h_1 (t) = \min_ {s \in [0,1]} \frac {G_m (t, s)} {G_m (1, s)}> 0, \hspace {1ex} \text {for all} \hspace {1ex} t \in \left (0,1 \right] $$ and $$ C_1 = \max_{t \in [0,1]} \left \{\max_{s \in [0,1]} \frac {G_m (t, s)} {G_m (1, s)} \right \} \geq 1, $$
from where the result follows.

\qed

\subsection {Case $ \gamma <0 $}

Next, we will give the expression of the Green's function of the problem
\begin{equation} \label{equation3_linear}
		u '' (t) -m ^ 2u (t) + \sigma (t) = 0, \quad 0 <t <1,
\end{equation}
coupled to the boundary conditions \eqref{condiciones}.

In this section we will omit most of the proofs because they are  analogous to those made in previous cases.

First, we will state a lemma that will be useful for the calculation of this function.

\begin{lemma} \cite[Appendix B]{Cabada_green}\label{green_homogeneo}
	Let's consider the problem
	$$ (P) \begin{cases}
		u '' (t) -m ^ 2u (t) + \sigma (t) = 0, \\
		u (0) = 0 = u (1).
	\end{cases} $$
	The Green's function associated with  problem $ (P) $ is given by the following expression
	\begin{equation} \label{function_green_ec3}
		G_m (t, s) = \begin{cases}
			\frac {\sinh (ms) \sinh (m (1-t))} {m \sinh m}, & 0 \le s \le t \le 1, \\
			\frac {\sinh (mt) \sinh (m (1-s))} {m \sinh m}, & 0 \le t <s \le 1.
		\end{cases}
	\end{equation}
\end{lemma}

\begin{theorem} \label{equation3_green_teorem}
	Let $ \lambda \neq \frac {m \sinh m} {\cosh m-1} $ and $ \sigma \in \mathcal {C} ([0,1]) $, then  problem (\ref{equation3_linear}),  (\ref{condiciones}) has a unique solution $ u \in \mathcal {C} ^ 2 ([0,1]) $, which is given by the expression $$ u (t) = \int_{0} ^ { 1} G_m (t, s) \sigma (s) ds, $$ where
	
	\begin{equation} \label{equation3_green}
		G_m (t, s) = \begin{cases}
			G ^ 1_m (t, s), & 0 \le s \le t \le 1, \\
			G ^ 2_m (t, s), & 0 \le t <s \le 1,
		\end{cases} \\
	\end{equation}
	with
	$$ \begin{aligned}
		G ^ 1_m (t, s) = & \frac {\sinh (ms) [\sinh (m-mt) (m \sinh m + \lambda (1- \cosh m)) - \lambda \sinh (mt)] } {m \sinh m (m \sinh m + \lambda (1- \cosh m))} \\& - \frac {\lambda \sinh (mt) (\sinh (m-ms) - \sinh m)} {m \sinh m (m \sinh m + \lambda (1- \cosh m))}
	\end{aligned} $$ and $$
	G ^ 2_m (t, s) = \frac {\sinh (mt) [\sinh (m-ms) (m \sinh m- \lambda \cosh m) - \lambda (\sinh (ms) - \sinh m )]} {m \sinh m (m \sinh m + \lambda (1- \cosh m))}. $$
\end{theorem}
\proof
In a similar way to Theorem \ref{equation2_green_teorem} and using  Lemma \ref{green_homogeneo}, we construct the Green's function taking into account that we can express the solution as
$$ u (t) = \int_{0} ^ {1} G_m (t, s) \sigma (s) ds + \lambda \left (\int_{0} ^ {1} u (s) ds \right) w (t), $$ where $ w (t) = \frac {\sinh mt} {\sinh m} $  is the unique solution of the following problem
$$ \begin{cases}
	w '' (t) -m ^ 2w (t) = 0, \quad 0 <t <1, \\
	w (0) = 0, w (1) = 1.
\end{cases} $$
\qed

We will now enunciate some properties of $ G_m $. The proofs are analogous to those that have been presented in the two previous cases.

\begin{lemma} \label{prop_green3}
	Let $ G_m $ be the Green's function associated with  problem (\ref{equation3_linear}),  (\ref{condiciones}), given by  expression (\ref{equation3_green}). Then for all $ \lambda \neq \frac {m \sinh m} {\cosh m-1} $, $m>0$, the following properties hold:
	\begin{enumerate}
		\item $ G_m (0, s) = G_m (t, 1) = G_m (t, 0) = 0 $, for all $ t, s \in [0,1]. $
		\item $ G_m (t, s) $ is continuous on $ [0,1] \times [0,1]. $
		\item $ G_m (1, s) = 0 $, for all $ s \in (0,1) $ if and only if $ \lambda = 0 $.
		\item $ (m \sinh m + \lambda (1- \cosh m)) G_m (1, s)> 0 $ for all $ s \in (0,1) $ and $m>0$.
		\item $G_m (t, s)> 0 $ for all $ t, s \in (0,1 )$ if and only if $ 0 \le \lambda <\frac {m \sinh m} {\cosh{m}-1}$, $m>0$. 
		\item $G_m (t, s) $ changes sign on $(0,1 ) \times (0,1)$  for all $ \lambda \not \in \left[0,\frac {m \sinh m} { \cosh m-1}\right)$, $m >0$.
	\end{enumerate}
\end{lemma}

\begin{lemma} \label{cotas_green3_enunc}
	Let $ m> 0 $, $ 0 \le \lambda <\frac {m \sinh m} {\cosh m-1} $ and $ G_m (t, s) $ be the Green's function of problem (\ref{equation3_linear}),  (\ref{condiciones}) given by  expression (\ref{equation3_green}). Then there are $ h_2 \in \mathcal {C} ([0,1]), \hspace {1ex} h_2> 0 $ in $ \left (0,1 \right] $ and $ C_2> 0 $ such that:
	\begin{equation} \label{cotas_green3}
		h_2 (t) G_m (1, s) \leq G_m (t, s) \leq C_2\,G_m (1, s), \qquad \text {for all}\; t, s \in [0,1]
	\end{equation}
\end{lemma}

\section{Nonlinear Problem}

This section is devoted to prove the existence of positive solutions on $(0,1)$ of problem \eqref{ecuacion_general}--\eqref{condiciones}. 
We will assume the following regularity condition for the nonlinear part of the equation:
\begin{flalign*}
	&(f)\qquad	f:[0,1]\times\left[0,\infty\right)\to\left[0,\infty\right)\hspace{2ex}\text{is a continuous function.}&
\end{flalign*}

As in the previous section, we will distinguish three different cases depending on the sign of the parameter $\gamma$. The results hold from the application of the Krasnosel'ski\u {\i}'s fixed point Theorem \ref{Krasnoselskii} to the operator $T_\gamma:\mathcal{C}([0,1])\to \mathcal{C}([0,1])$ defined as
\begin{equation}\label{operador_teorema}
	T_\gamma u(t):=\int_{0}^{1}G_\gamma(t,s)f(s,u(s))ds, \quad t \in [0,1].
\end{equation}

Here,  $G_\gamma$ corresponds to function $G$, given by \eqref{ecuacion1_green}, if $\gamma=0$, $G_m$ given by \eqref{equation2_green} if $\gamma=m^2>0$, and $G_m$ given by \eqref{equation3_green} if $\gamma=-m^2<0$.

As we have proved in previous section, we know that the fixed points of  operator $ T_\gamma $ coincide with the solutions of  problem (\ref{ecuacion_general})--(\ref{condiciones}). 

To apply Theorem \ref{Krasnoselskii}, we define $ X = (\mathcal{C}( [0,1]), || \cdot ||), $ the Banach space endowed with the supremum norm.

Now, we denote
$$f_0=\lim\limits_{u\to 0^+}\left\{\min_{t\in[\frac{1}{2},1]} \frac{f(t,u)}{u}\right\}, \qquad\text{}\qquad f_\infty=\lim\limits_{u\to +\infty}\left\{\min_{t\in[\frac{1}{2},1]} \frac{f(t,u)}{u}\right\}$$
and
$$f^0=\lim\limits_{u\to 0^+}\left\{\max_{t\in[0,1]} \frac{f(t,u)}{u}\right\} \qquad\text{and}\qquad f^\infty=\lim\limits_{u\to +\infty}\left\{\max_{t\in[0,1]} \frac{f(t,u)}{u}\right\}.$$

In the sequel, we introduce the cone $K_\gamma\subset X$, depending on the sign of the real parameter $\gamma$.
	
If $\gamma=0$:
	\begin{equation}\label{cone1}
		K_\gamma=\left\{u\in X;\hspace{1ex}u(t)\geq0,\;t\in[0,1], \hspace{1ex} u(t)\geq t\frac{\lambda}{2}||u||, \;t\in\left[0,1\right] \right\}.
	\end{equation}
	
	If $\gamma>0$:
			\begin{equation} \label{cone2}
		K_\gamma = \left \{u \in X; \hspace {1ex} u (t) \geq0, \hspace {1ex} u (t) \geq \frac {h_1 (t)} {C_1} || u || , \; t \in \left [0,1 \right] \right \}.
	\end{equation}

	If $\gamma<0$:
\begin{equation} \label{cone3}
	K_\gamma = \left \{u \in X; \hspace {1ex} u (t) \geq0, \hspace {1ex} u (t) \geq \frac {h_2 (t)} {C_2} || u || , \; t \in \left [0,1 \right] \right \}.
\end{equation}

Here $h_1$ and $C_1$ are given in Lemma \ref{cotas_green2_enunc}, and $h_2$ and $C_2$ are given in Lemma \ref{cotas_green3_enunc}.

	So, we arrive at the following existence result.
	
\begin{figure} [h]
		\centering
		\includegraphics[scale = 0.75]{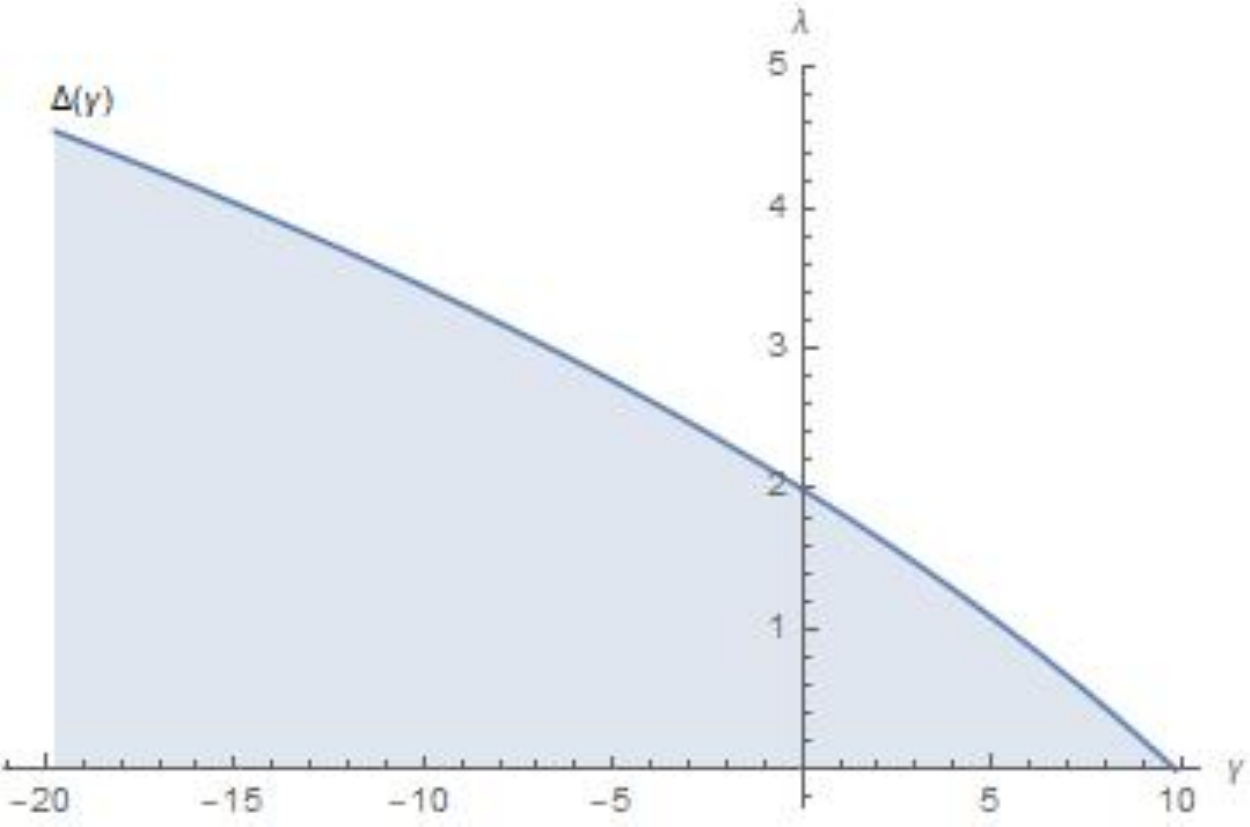}
		\caption {Graphic of $ \Delta (\gamma) $ \label{fig-delta}}
	\end{figure}

		\begin{theorem} \label{general_solution_equation}
		Let us consider  problem (\ref{ecuacion_general})--(\ref{condiciones}), and let $ \Delta: (- \infty, \pi ^ 2) \to \mathbb {R} $ be the function defined as (see Figure \ref{fig-delta}):
		\begin{equation} \label{delta_function}
			\Delta (\gamma) = \begin{cases}
				\frac {\sqrt {- \gamma} \sinh (\sqrt {- \gamma})} {\cosh (\sqrt {- \gamma}) - 1}, & \gamma <0, \\
				2, & \gamma = 0, \\
				\frac {\sqrt {\gamma} \sin (\sqrt {\gamma})} {1- \cos (\sqrt {\gamma})}, & \gamma >0.
			\end{cases}
		\end{equation}
		Suppose further that $ (f) $ holds and one of the two following conditions is fulfilled:
		\begin{enumerate}
			\item [(i)] (sublinear case) $ f_0 = \infty $ and $ f ^ \infty = 0. $
			\item [(ii)] (superlinear case) $ f ^ 0 = 0 $ and $ f_ \infty = \infty. $
		\end{enumerate}
		So, for all $ \gamma <\pi ^ 2$ and $ 0 <\lambda <\Delta (\gamma) $ there is a positive solution of problem (\ref{ecuacion_general})--(\ref{condiciones}), $u\in K_\gamma$.
	\end{theorem}
	
\proof  Consider, in a first moment, the case $\gamma=m^2>0$.

		Let's first see that $ T: K_\gamma \to K_\gamma $ is a compact operator. \\
		Since $ G_m $ and $ f $ are continuous and non-negative in their domain of definition, then $ T_\gamma u \in \mathcal {C} ([0,1]) $ and $ T_\gamma u (t) \geq0 $, for all $ t \in [0,1]. $ \\
		Let $ u \in K_\gamma $, using the properties stated in  lemmas \ref{prop_green2} and \ref{cotas_green2_enunc}, we have that, for all $ t \in [0,1] $
		$$ \begin{aligned}
		T_\gamma u (t) & = \int_{0} ^ {1} G_m (t, s) f (s, u (s)) ds \geq \int_{0} ^ {1} h_1 (t) G (1, s) f (s, u (s)) ds \\& \geq \frac {h_1 (t)} {C_1} \int_{0} ^ {1} \left \{\max_{t \in [0, 1]} G (t, s) \right \} f (s, u (s)) ds \\ & \geq \frac {h_1 (t)} {C_1} \max_{t \in [0,1]} \left \{\int_{0} ^ {1} G_m (t, s) f (s, u (s)) ds \right \} \\& = \frac {h_1 (t)} {C_1} || T_\gamma u | |.
				\end{aligned} $$
				Furthermore, the continuity of the functions $ G_m $ and $ f $ guarantees the continuity of the operator $ T: K_\gamma \to K_\gamma $. So $ T (K_\gamma) \subset K_\gamma $. \\
				Let us now verify that the image by $ T_\gamma $ of a bounded set is relatively compact. To this end, we will use the Arzel\`{a}-Ascoli Theorem. \\
				Let $ \Omega \subset K_\gamma $ bounded, that is, there exists $ M \in \mathbb {R}, \hspace {1ex} M> 0 $ such that $ || u || \leq M $, for all $ u \in \Omega. $ Let's define $$ L = \max_{0 \leq t \leq1,0 \leq u \leq M} | f (t, u) |. $$
				Then for all $ u \in \Omega $ and $ t \in [0,1] $, we have that
				$$ \begin{aligned}
				| T_\gamma u (t) | = & \left | \int_{0} ^ {1} G_m (t, s) f (s, u (s)) ds \right | \leq L \int_{0} ^ {1 } G_m (t, s) ds \leq L \int_{0} ^ {1} C_1G_m (1, s) ds \\
				= & LC_1 \frac {\lambda (2-2 \cos m-m \sin m)} {m ^ 2 (m \sin m- \lambda (1- \cos m))}: = N.
				\end{aligned} $$
				As a consequence, $$ || T_\gamma u || \leq N, $$ so $ T _\gamma(\Omega) $ is bounded. \\
				On the other hand, for each $ u \in \Omega $ and $ t \in [0,1] $ we have
				$$ \begin{aligned}
				| (T_\gamma u) '(t) | = & \left | \int_{0} ^ {1} \frac {\partial G_m} {\partial t} (t, s) f (s, u (s)) ds \right | \leq
				\int_{0} ^ {1} \left | \frac {\partial G_m} {\partial t} (t, s) \right || f (s, u (s)) | ds \\\leq & L \int_{0} ^ {1} \left | \frac {\partial G_m} {\partial t} (t, s) \right | ds =: N '.
				\end{aligned} $$
				The regularity of the Green's function allows us to ensure that $ N ' \in \R $, $N'>0$. Thus, for all $ t_1, t_2 \in [0,1], \hspace {1ex} t_1 <t_2 $, the following inequality is fulfilled
				$$ | (T_\gamma u) (t_2)-(T_\gamma u) (t_1) | = \left | \int_{t_1} ^ {t_2} (T_\gamma u) '(s) ds \right | \leq \int_{t_1} ^ { t_2} | (T_\gamma u) '(s) | ds \leq N' (t_2-t_1), $$
				so $ T _\gamma(\Omega) $ is an equicontinuous set in $ X $. \\
				Thus, by virtue of the Arzel\`{a}-Ascoli Theorem, we deduce that the set $ \overline {T _\gamma(\Omega)} $ is relatively compact, i.e., $ T: K_\gamma \to K_\gamma $ is a compact operator. \\
			
		Suppose that the first situation holds:
		\begin{enumerate}
		\item [(i)] (sublinear case) $ f_0 = \infty $ and $ f ^ \infty = 0. $
		\end{enumerate}
		Since $ f_0 = \infty $, there exists a constant $ \rho_1 $ such that $ f (t, u) \geq \delta_1 u $, for all $t \in [1/2,1]$ and $ 0 \leq u \leq \rho_1 $, where $ \delta_1 $ is such that $$ \frac {\delta_1} {C_1} \max_{t \in [0,1]} \left \{\int_{\frac {1} {2}} ^ {1} h_1 (s) G_m (t, s) ds \right \} \geq 1. $$
			Let $ u \in K_\gamma $ with $ || u || = \rho_1 $, then from the previous expression we deduce
			$$ \begin{aligned}
			|| T_\gamma u || & = \max_{t \in [0,1]} \left \{\int_{0} ^ {1} G_m (t, s) f (s, u (s)) ds \right \} \geq \max_{t \in [0,1]} \left \{\int_{\frac {1} {2}} ^ {1} G_m (t, s) f (s, u (s) ) ds \right \} \\& \geq
					\max_{t \in [0,1]} \left \{\int_{\frac {1} {2}} ^ {1} G_m (t, s) \delta_1u (s) ds \right \} \geq \delta_1 \max_{t \in [0,1]} \left \{\int_{\frac {1} {2}} ^ {1} G_m (t, s) \frac {h_1 (s)} {C_1 } || u || ds \right \} \\& = || u || \frac {\delta_1} {C_1} \max_{t \in [0,1]} \left \{\int_{\frac {1} {2}} ^ {1} h_1 (s) G_m (t, s) \right \} \geq || u ||.
								\end{aligned} $$
								On the other hand, the continuity of $ f $ in the second variable allows us to define the function
								$$ \widetilde {f} (t, u) = \max_{z \in [0, u]} \{f (t, z) \}, \quad t \in [0,1], \hspace { 1ex} u \in \mathbb {R}, $$
									which is monotone nondecreasing on $ \left [0, \infty \right) $ for every $ t \in [0,1] $. Now, since $ f ^ \infty = 0 $, it follows (see \cite {Wang}) that
								\begin{equation}
									\label{e-ftilde}
									\lim \limits_ {u \to \infty} \left \{\max_{t \in [0,1]} \frac {\widetilde {f} (t, u)} {u} \right \} = 0. 
								\end{equation} 
										Let's now take $ \delta_2> 0 $ be such that $$ \delta_2 \max_{t \in [0,1]} \left \{\int_{0} ^ {1} G_m (t, s) ds \right \} \leq1. $$ 
										
										From \eqref{e-ftilde}, we know that there exists $ \rho_2 \in \mathbb {R}, \hspace {1ex} \rho_2> \rho_1> 0 $ such that $ \widetilde {f} (t, u) \leq \delta_2u $, for all $t \in [0,1]$ and $ u \geq \rho_2. $ \\
											Let $ u \in K_\gamma $ be such that $ || u || = \rho_2 $, then using the definition of $ \widetilde {f} $ and the above inequality we have
											$$ \begin{aligned}
											|| T_\gamma u || & = \max_{t \in [0,1]} \left \{\int_{0} ^ {1} G_m (t, s) f (s, u (s)) \right \} \leq \max_{t \in [0,1]} \left \{\int_{0} ^ {1} G_m (t, s) \widetilde {f} (s, u (s)) \right \} \\& \leq \max_{t \in [0,1]} \left \{\int_{0} ^ {1} G_m (t, s) \widetilde {f} (s, || u || ) \right \} \leq \max_{t \in [0,1]} \left \{\int_{0} ^ {1} G_m (t, s) \delta_2 || u || \right \} \\& \leq || u || \delta_2 \max_{t \in [0,1]} \left \{\int_{0} ^ {1} G_m (t, s) \right \} \leq || u ||.
														\end{aligned} $$
														Finally, the first part of  Theorem \ref{Krasnoselskii}  implies that there is at least one positive solution of problem (\ref{ecuacion_general})--(\ref{condiciones}), with $\gamma=m^2>0$, $ u \in K_\gamma $, such what
														$ \rho_1 \leq || u || \leq \rho_2. $ \\
			
			Let us now consider the second case:
			\begin{enumerate}
			\item [(ii)] (superlinear case) $ f ^ 0 = 0 $ and $ f_ \infty = \infty. $
			\end{enumerate}
			Let $ \delta_2 \in \left (0,1 \right] $, chosen in the same way as in the previous case. Since $ f ^ 0 = 0 $, there exists a constant $ r_1> 0 $ such that $ f ( t, u) \leq \delta_2u $ for all $t \in [0,1]$ and $ 0 \leq u \leq \delta_2 $. \\
			Let's choose $ u \in K_\gamma $ be such that $ || u || = r_1 $, then
			$$ \begin{aligned}
			|| T_\gamma u || & = \max_{t \in [0,1]} \left \{\int_{0} ^ {1} G_m (t, s) f (s, u (s)) \right \} \leq \max_{t \in [0,1]} \left \{\int_{0} ^ {1} G_m (t, s) \delta_2 u (s) \right \} \leq \\& \leq || u || \delta_2 \max_{t \in [0,1]} \left \{\int_{0} ^ {1} G_m (t, s) \right \} \leq || u || .
					\end{aligned} $$
					Let's take $ \delta_3> 0 $ so that
					$$ \frac {\delta_3} {C_1} \max_{t \in [0,1]} \left \{\int_{\frac {1} {2}} ^ {1} h_1 (s) G_m (t , s) ds \right \} \geq1. $$
						Since $ f_ \infty = \infty $, there exists $ r_2> r_1> 0 $ with $ C_1r_2> \underline {M} r_1 $ (where $ \underline {M} = \min_ {t \in [\frac {1 } {2}, 1]} h_1 (s) $), so that $ f (t, u) \geq \delta_3 u $ for all $t \in [1/2,1]$ and $ u \geq r_2. $
						Let $ u \in K_\gamma $ be such that $ || u || = r_2 \frac {C_1} {\underline {M}} $, then, by the definition of $ K_\gamma $ it is satisfied that $ u (t) \geq r_2 $ for all $ t \in \left [1/2, 1 \right] $. \\
						From this, we deduce the following inequalities
						$$ \begin{aligned}
						|| T_\gamma u || & = \max_{t \in [0,1]} \left \{\int_{0} ^ {1} G_m (t, s) f (s, u (s)) ds \right \} \geq
							\max_{t \in [0,1]} \left \{\int_{\frac {1} {2}} ^ {1} G_m (t, s) f (s, u (s)) ds \right \} \\& \geq \max_{t \in [0,1]} \left \{\int_{\frac {1} {2}} ^ {1} G_m (t, s) \delta_3u (s) ds \right \} \geq \delta_3 \max_{t \in [0,1]} \left \{\int_{\frac {1} {2}} ^ {1} G_m (t, s) \frac { h_1 (s)} {C_1} || u || ds \right \} \\& = || u || \frac {\delta_3} {C_1} \max_{t \in [0,1]} \left \{\int_{\frac {1} {2}} ^ {1} h_1 (s) G_m (t, s) ds \right \} \geq || u ||.
											\end{aligned} $$
											Finally, applying part 2 of  Theorem \ref{Krasnoselskii}, we conclude that problem (\ref{ecuacion_general})--(\ref{condiciones}), with $\gamma=m^2>0$, has at least one positive solution $ u \in K_\gamma $ such that
											$$ r_1 \leq || u || \leq \frac {C_1} {\underline {M}} r_2. $$

											The cases $\gamma=0$ and $\gamma=-m^2<0$ can be proved in a similar way.
											\qed

\begin{remark}
	\label{r-a-b}
It should be noted that, since $ G (0, s) = 0 $, in order to ensure the existence of $ r_2 $ in the case $ (ii) $ of Theorem \ref{general_solution_equation}, we have reduced the interval considered in the definition of the limits $ f_0 $ and $ f_ \infty $ from $ [0,1] $ to $ \left [1/2, 1 \right ] $. \\
In fact, it is enough to take $ a, b \in \mathbb {R} $ such that $ 0 <a <b \le 1 $ and redefine the limits as
$$ f_0 = \lim \limits_ {u \to 0 ^ +} \left \{\min_ {t \in [a, b]} \frac {f (t, u)} {u} \right \} \qquad \text {and} \qquad f_ \infty = \lim \limits_ {u \to + \infty} \left \{\min_ {t \in [a, b]} \frac {f (t, u)} {u} \right \}. $$
		In this case, it is easy to check that Theorem \ref{general_solution_equation} remains true for solutions defined in cone $K_\gamma$, given by expressions \eqref{cone1},  \eqref{cone2} and \eqref{cone3}, by replacing on their definitions $[1/2,1]$ by $[a,b]$.
			\end{remark}

			\begin{remark}
			Let us now consider  problem (\ref{ecuacion_general}) coupled to the boundary conditions
			\begin{equation} \label{conditions2}
			u (0) = \lambda \int_{0}^{1} u (s) ds, \qquad u (1) = 0.
		\end{equation}
		It is easy to check that $ v (t) :=u (1-t)$, $t \in [0,1]$, also satisfies (\ref{ecuacion_general}) together with
		$$ v (0) = u (1) = \lambda \int_{0} ^ {1} u (s) ds = \lambda \int_{0} ^ {1} v (s) ds, \qquad v ( 1) = u (0) = 0. $$
		\end{remark}
	
	So, if we denote
	$$f'_0=\lim\limits_{u\to 0^+}\left\{\min_{t\in[0,\frac{1}{2}]} \frac{f(t,u)}{u}\right\} \qquad\text{and}\qquad f'_\infty=\lim\limits_{u\to +\infty}\left\{\min_{t\in[0,\frac{1}{2}]} \frac{f(t,u)}{u}\right\}$$
	and, if $\gamma=0$:
	\begin{equation*}
		K'_\gamma=\left\{u\in X;\hspace{1ex}u(t)\geq0,\;t\in[0,1], \hspace{1ex} u(t)\geq (1-t)\frac{\lambda}{2}||u||, \;t\in\left[0,1\right] \right\}.
	\end{equation*}
	
	If $\gamma>0$:
	\begin{equation*}
		K'_\gamma = \left \{u \in X; \hspace {1ex} u (t) \geq0, \hspace {1ex} u (t) \geq \frac {h_1 (1-t)} {C_1} || u || , \; t \in \left [0,1 \right] \right \}.
	\end{equation*}
	
	If $\gamma<0$:
	\begin{equation*}
		K'_\gamma = \left \{u \in X; \hspace {1ex} u (t) \geq0, \hspace {1ex} u (t) \geq \frac {h_2 (1-t)} {C_2} || u || , \; t \in \left [0,1 \right] \right \}.
	\end{equation*}

		Thus, the following result is then obtained.

		\begin{corollary}
		Consider  problem (\ref{ecuacion_general}),  (\ref{conditions2}) and let $ \Delta: (- \infty, \pi ^ 2) \to \mathbb {R} $ be the function defined in (\ref{delta_function}).
		Suppose further that $ (f) $ holds and one of the following conditions is satisfied:
		\begin{enumerate}
		\item [(i)] (sublinear case) $ f'_0 = \infty $ and $ f ^ \infty = 0. $
\item [(ii)] (superlinear case) $ f ^ 0 = 0 $ and $ f'_ \infty = \infty. $
		\end{enumerate}
				So, for all $ \gamma <\pi ^ 2$ and $ 0 <\lambda <\Delta (\gamma) $ there is a positive solution of problem (\ref{ecuacion_general}),  (\ref{conditions2}), $u\in K'_\gamma$.
	\end{corollary}
	
	Moreover, analogous considerations to Remark \ref{r-a-b} remain valid for problem (\ref{ecuacion_general}),  (\ref{conditions2}).
	\subsection {Examples}
	
	\begin{example}
	Consider  problem (\ref{ecuacion_general})--(\ref{condiciones}) with
	$$ f (t, x) = \sqrt [5] {x ^ 3 +x} + \log (3\,t + x). $$
	It is easy to verify that, for $ u> 0 $
	$$ \min_ {t \in [\frac {1} {2}, 1]} \frac {f (t, u)} {u} = \frac {\sqrt [5] {u ^ 3 + u} + \log (\frac {3} {2} + u)} {u} \; \text {and} \; \max_ {t \in [0,1]} \frac {f (t, u)} {u} = \frac {\sqrt [5] {u ^ 3 + u} + \log (3 + u)} {u}. $$
	Taking limits, it is easy to check that $ f_0 = \infty $ and $ f ^ \infty = 0 $, i.e., we are in the sublinear case. Since $ f $ fulfills $ (f) $, we are in the hypotheses of the case (i) Theorem \ref{general_solution_equation}, so we can ensure the existence of a positive solution for problem (\ref{ecuacion_general})-( \ref{condiciones}).
	\end{example}
	
	\begin{example}
	Consider again the problem (\ref{ecuacion_general})--(\ref{condiciones}) with
	$$ f (t, x) = t\, x^ 3 + e ^ {t\,x}-1.$$
	In this case, we have for all $ u> 0 $
	$$ \min_ {t \in [\frac {1} {2}, 1]} \frac {f (t, u)} {u} = \frac {u ^ 2} {2} + \frac {e ^ {\frac {u ^ 2} {2}} - 1} {u} \; \mbox{and} \, \max_ {t \in [0,1]} \frac {f (t, u)} {u} = u ^ 2 + \frac {e ^ {u ^ 2} -1} {u}. $$
	Again, taking limits we see that $ f ^ 0 = 0 $ and $ f_ \infty = \infty $, that is, we are in the superlinear case. As $ f $ verifies $ (f) $, the second part of  Theorem \ref{general_solution_equation} guarantees the existence of a positive solution of  problem (\ref{ecuacion_general})--(\ref{condiciones}).
	\end{example}

	

	\end{document}